\documentclass[11pt]{amsart}

\usepackage{amsmath, amsfonts, amssymb}

\usepackage{hyperref}

\newcommand{\C}{\mathbb{C}}
\newcommand{\T}{\mathbb{T}}
\newcommand{\D}{\mathbb{D}}

\DeclareMathOperator{\VMO}{VMO}

\DeclareMathOperator{\IDA}{IDA}
\DeclareMathOperator{\VDA}{VDA}
\DeclareMathOperator{\ind}{ind}
\DeclareMathOperator{\wind}{wind}
\DeclareMathOperator{\cl}{cl}

\title{On the product formula for Toeplitz and related operators}

\author{Jani A. Virtanen}

\address{Department of Mathematics, University of Reading, England}

\email{j.a.virtanen@reading.ac.uk}

\dedicatory{In memory of Harold Widom}



\begin{document}

\begin{abstract} In this note known formulas for the product of Toeplitz operators are revisited in the context of their applications to the study of Fredholmness, boundedness of Toeplitz products, and the Berezin-Toeplitz quantization. A few open problems are also mentioned.
\end{abstract}

\maketitle

\section{Introduction}
Given two bounded Toeplitz operators $T_f$ and $T_g$ on the Hardy space $H^2$, their product can be written as
\begin{equation}\label{e:1}
	T_f T_g = T_{fg} - H_f H_{\tilde g},
\end{equation}
where $H_f$ and $H_{\tilde g}$ are Hankel operators acting on $H^2$. As stated in~\cite{BBE}, this identity was established by Widom~\cite{Widom}, while it had been known and used for a long time in other forms, such as 
$$
	PfPgP=PfgP-PfQgP, 
$$
where $P$ is the orthogonal projection of $L^2$ onto $H^2$ and $Q = I-P$. 
What resulted from Widom's use of this identity was a very ingenious way of dealing with the asymptotics of block Toeplitz determinants in~\cite{Widom}, now known as the Szeg\"o-Widom asymptotics, 
via operator theoretic methods and Schatten class properties of Hankel operators. Paper \cite{BBE} is embarking on this topic.

Going back to the identity in \eqref{e:1} and its original intent to show that certain Toeplitz operators are Fredholm, I will discuss extensions of this formula in the context of other function spaces, such as Bergman and Fock spaces, and show how it leads to interesting questions about the properties of Hankel operators. What we lack in these other function spaces, however, are effective matrix representations of Toeplitz and Hankel operators, which creates an obstacle to obtaining Widom type identities for the products of truncated Toeplitz matrices.

For simplicity, we limit the discussion to function spaces defined over domains in $\C$, except for Section~\ref{quantization}, and note that the generalizations to the $n$-dimensional setting can be easily found in the literature. 

\section{Preliminaries}

For $0<p<\infty$, $\Omega\subset\C$, and $\mu$ a positive measure on $\Omega$, denote by $L^p(\Omega, d\mu)$ the space of all complex measurable functions $f$ on $\Omega$ for which
$$
	\|f\|_p = \left( \int_\Omega |f|^p \, d\mu \right)^{1/p} < \infty.
$$
For a complex measurable function $f$ on $\Omega$, define $\|f\|_\infty$ to be the essential supremum of $|f|$ and denote by $L^\infty(\Omega, d\mu)$ all $f$ for which $\|f\|_\infty<\infty$. The set of all analytic functions in an open set $\Omega$ is denoted by $H(\Omega)$.

In terms of domains $\Omega$, the usual three model cases consist of the unit circle $\T$, the unit disk $\D$, and the complex plane $\C$. When $\Omega = \T$, we write $L^p(\T)$ for $L^p(\T, d\theta)$ and define the Hardy space $H^p$ by
$$
	H^p = \{ f\in L^p(\T) : f_k = 0\ {\rm for}\ k<0\}.
$$
Let $dA=dxdy$ be the usual area measure on $\C$. We write $L^p(\D)$ for $L^p(\D, dA)$ and define the Bergman space $A^p$ by
$$
	A^p = H(\D) \cap L^p(\D).
$$
When $\Omega=\C$, define the Fock space $F^p$ by
$$
	F^p = H(\C) \cap L^p(\C, e^{-\frac{p}2|z|^2}dA).
$$

Let $X^2(\Omega) \in \{ H^2, A^2, F^2\}$. Then $X^2(\Omega)$ is a Hilbert space and the orthogonal projection of $L^2(\Omega)$ onto $X^2(\Omega)$ is denoted by $P$. We write $Q = I-P$ for the complementary projection. Given a bounded function $f$ on $\Omega$, the Toeplitz operator $T_f : X^p(\Omega) \to X^p(\Omega)$ with symbol $f$ is defined by
$$
	T_fg = P(fg).
$$
When $1<p<\infty$, since $P$ extends to a bounded projection on $L^p(\Omega)$, $T_f$ is clearly bounded on $X^p(\Omega)$ if $f$ is bounded.

Defining Hankel operators is less straightforward. Indeed, the Hankel operators that appear in~\eqref{e:1} act on the Hardy space while the Hankel operators on Bergman spaces $A^p$ or Fock spaces $F^p$ map into the corresponding $L^p(\Omega)$. More precisely, define the flip operator $J : L^p(\T) \to L^p(\T)$ by 
$$
	Jf(t) = \bar t f(\bar t)
$$
for $t\in \T$. For a bounded symbol $f$, the Hankel operator $H_f$ is defined on $H^p$ by
$$
	H_fg = PM_f QJf,
$$
where $M_f$ is the multiplication operator. When $\Omega\in \{\D, \C\}$, we define the Hankel operator $H_f : X^p(\Omega) \to L^p(\Omega)$ by
$$
	H_fg = Q(fg).
$$
Again, it is easy to see that the Hankel operator $H_f$ is bounded in all the three cases if $1<p<\infty$ and $f$ is bounded.

\section{Fredholm properties of Toeplitz opeators}
In this section the Fredholm properties of Toeplitz operators acting on Hardy, Bergman and Fock spaces are considered using~\eqref{e:1} and its generalizations. Recall that an operator $A$ on a Banach space is said to be Fredholm if $\ker A$ and $X/A(X)$ are both finite dimensional, in which case the index $\ind A$ is defined by
$$
	\ind A = \dim \ker A - \dim X/A(X).
$$
Equivalently, $A$ is Fredholm if and only if $A+K(X)$ is invertible in the Calkin algebra $B(X)/K(X)$, where $B(X)$ and $K(X)$ denote the sets of all bounded and compact operators on $X$, respectively. The essential spectrum of $A$ is defined by
$$
	\sigma_{\rm ess} (A) = \{\lambda \in \C : A-\lambda\ \textrm{is not Fredholm}\}.
$$	

\subsection{The Hardy space case}
Let $f, g$ be bounded on $\T$ and write $\tilde f(t) = f(\bar t)$ for $t\in \T$. Then
\begin{equation}\label{e:1 proof}
\begin{split}
	T_{fg} &= PM_{fg}P = PM_fM_gP = PM_fPM_bP + PM_fQM_gP\\
	&=PM_fP^2 M_g P + PM_fQJ^2 QM_gP,
\end{split}
\end{equation}
which is \eqref{e:1}.

Suppose now that $f$ is continuous and has no zeros on $\T$. Then $g=1/f$ is also continuous and has no zeros. By~\eqref{e:1}, since $H_f$ is known to be compact,
$$
	T_f T_g = I + H_f H_{\tilde g} = I + K
$$
for some compact operator $K$. Similarly, $T_gT_f - I$ is compact, and hence $T_f$ is Fredholm. In situations when Hankel operators are compact, the identity in~\eqref{e:1} is tailor-made for proving that Toeplitz operators are Fredholm. In other words, whenever the Hankel operators are compact, the corresponding Toeplitz operators commute modulo compact operators. A similar approach also applies to symbols in the Douglas algebra $C+\overline{H^\infty}$ but the use of~\eqref{e:1} is no longer as effective with more general classes of symbols.

Let $f\in L^\infty(\T)^{N\times N}$  and consider the block Toeplitz operator $T_f$ on $H^p_N = \{ (f_1,\ldots, f_N)^\top : f_j \in H^p\}$. Suppose that $f\in (C+\overline{H^\infty})^{N\times N}$ and $\det f$ is invertible in $C + \overline{H^\infty}$. Choose $h \in (\mathcal{R} + \overline{H^\infty})^{N\times N}$, where $\mathcal{R}$ is the set of all rational functions, sufficiently close to $f$ in the norm of $L^\infty_{N\times N}(\T)$. Then
$$
	\ind T_f = \ind T_h \quad{\rm and}\quad \ind T_{\det f} = \ind T_{\det h}.
$$
Since $H_h$ has finite rank, \eqref{e:1} implies that the entries of $T_h$ commute modulo finite-rank operators, and hence $\ind T_h = \ind T_{\det h}$ (see Theorem~1.15 of~\cite{BS2006}), which reduces the index computation to that of the scalar-valued symbols. For more general symbols, such as piecewise continuous symbols, no such reductions are possible.

Although the projection $P$ is unbounded on $L^1(\T)$, a Fredholm theory for Toeplitz operators on $H^1$ can still be developed. In particular, when $f$ is a continuous function of logarithmic vanishing mean oscillation, the Hankel operator $H_f$ is compact on $H^1$, so~\eqref{e:1} is readily available, and the Fredholm properties can be described as in the reflexive case $1<p<\infty$ (see~\cite{V2006}). However, as recently observed, there are continuous symbols $f$ that generate bounded Toeplitz operators on $H^1$ and for which $H_f$ is not compact (see~\cite{HV2022}). This makes the study of Fredholmness of $T_f$ with such continuous symbols considerably more difficult in $H^1$ because~\eqref{e:1} no longer produces a desired conclusion.

There are many other aspects of Toeplitz operators on the Hardy space whose proofs benefit from~\eqref{e:1}, such as invertibility with analytic symbols, the applicability of local principles, the study of Toeplitz algebras and Fisher-Hartwig symbols, but we refrain from further details (all of which can be found in~\cite{BS2006}) and keep our focus only on the Fredholm properties in this section.

\subsection{The Bergman space case}
As mentioned above, for $1<p<\infty$ and $f\in L^\infty(\D)$, the Hankel operator $H_f$ is defined by $H_fg = Q(fg)$ for $g\in A^p$, and so it maps into $L^p(\D)$ instead of $A^p$. However, we can still obtain formulas similar to~\eqref{e:1} as follows. For two bounded functions $f,g$ on $\Omega$, using the inner product in $A^2$, it is easy to see that
\begin{equation}\label{e:3}
	T_f T_g = T_{fg} - H^*_{\bar f} H_g
\end{equation}
when $p=2$, which shows that
$$
	T_{|f|^2} - T_{\bar f}T_f = H^*_f H_f,
$$
and hence compactness of $H_f$ is equivalent to compactness of the semi-self-commutator $T_{|f|^2} - T_{\bar f}T_f$. In addition, the formulas
\begin{align}\label{e:4}
	T_f T_g &= PM_f PM_g = PM_f(I-Q)M_g = T_{fg} - PM_fH_g\\
	\label{e:5}
	&=I - P(I-M_{fg}) - PM_fH_g = I - T_{1-fg} - PM_f H_g
\end{align}
are useful. For example, in~\cite{McDonald}, the identity in~\eqref{e:4} was used to show that the Toeplitz operator $T_f$ with $f\in C(\overline{\D})$ is Fredholm on $A^2$ if and only if $f$ has no zeros on the boundary. A similar approach, using~\eqref{e:5}, can be used to treat symbols in the Douglas algebra $C(\overline{\D}) + H^\infty$ and symbols of vanishing mean oscillation.

Let $f\in L^\infty(\D)^{N\times N}$ and consider the block Toeplitz operator $T_f$ on $A^p_N = \{ (f_1, \ldots, f_N)^\top : f_j \in A^p\}$. Fredholmness of block Toeplitz operators with symbols in the Douglas algebra $(C(\overline{\D}) + H^\infty(\D))^{N\times N}$ can be handled as in the Hardy space case but now with the identities in~\eqref{e:4} and~\eqref{e:5}. However, the index formula for these symbols cannot be derived as easily as in the Hardy space case because the formula $\ind A = \ind \det A$, which holds for operator matrices $A$ whose entries commute modulo trace class operators, fails to reach all of $C(\overline{\D}) + H^\infty(\D)$ via~\eqref{e:4}. For an alternate approach to the computation of the index of block Toeplitz operators $T_f$ on the Bergman spaces on the unit ball, see~\cite{BP2013}. Similar comments can be made about symbols in $(L^\infty(\D) \cap \VMO)^{N\times N}$, where $\VMO$ is the space of functions of vanishing mean oscillation, and in particular the approach in~\cite{BP2013} should produce an index formula for this symbol class, too. We return to this topic in the next section when dealing with Toeplitz operators in the Fock space setting, where an analogous problem still remains open. 

\subsection{The Fock space setting} As in the Bergman space case, for Toeplitz operators on $F^2$, we again have
\begin{equation}\label{e:6}
	T_fT_g = T_{fg} - H^*_{\bar f} H_g.
\end{equation}
To my knowledge, this identity was first used to describe the Fredholm properties of $T_f$ on the Fock space in~\cite{Stroethoff}. It was shown that, when $f\in L^\infty(\C)$ and $H_f$ is compact, we have
\begin{equation}\label{e:essential spectrum}
	\sigma_{\rm ess} (T_f) = \bigcap_{r>0} \cl \tilde f (\C\setminus \D_r),
\end{equation}
where $\cl E$ stands for the closure of $E$ in $\C$, $\D_r = \{|z|<r\}$, and $\tilde f$ is the Berezin transform of $f$ defined by
\begin{equation}\label{e:Berezin}
	\tilde f(z) = \frac1{2\pi} \int_{\C} f(w) e^{-\frac12|z-w|^2} dA(w)
\end{equation}
for $z\in\C$. In the proof of~\eqref{e:essential spectrum}, identity~\eqref{e:6} comes into play as follows. Suppose that $\xi\notin \cl f(\C\setminus \D_r)$ for some $r>0$. To show that $T_{f-\xi}$ is Fredholm, define
$$
	g(z) = \begin{cases}
	(f(z) - \xi)^{-1} & {\rm if}\ z\in \C\setminus \D_r,\\
	1 & {\rm if}\ z\in \D_r.
	\end{cases}
$$
Then $g\in L^\infty(\C)$, and an application of~\eqref{e:6} shows that 
$$
	T_{g}T_{f-\xi} = I - H^*_{\bar g} H_f - T_{(f-\xi-1)\chi_{\D_r}}.
$$
Notice that $(f-\xi-1)\chi_{\D_r}$ has compact support and hence $T_{(f-\xi-1)\chi_{\D_r}}$ is compact. Since $H_f$ is compact, it follows that $T_{f-\xi} + K(F^2)$ is left-invertible in $B(F^2)/K(F^2)$. That $T_{f-\xi} + K(F^2)$ is also right-invertible follows from $T_{f-\xi} = T^*_{\bar f - \bar\xi}$ and the fact that $H_{\bar f}$ is compact whenever $H_f$ is compact (see, e.g.,~\cite{BC1987} or~\cite{HuV2022-2}). Therefore, $T_{f-\xi} = T_f - \xi$ is Fredholm, that is, $\xi \notin \sigma_{\rm ess}(T_f)$, and so $\sigma_{\rm ess} (T_f) \subset \cl f(\C\setminus D_r)$ for all $r>0$. Further, since $T_{f-\tilde f}$ is known to be compact, $\sigma_{\rm ess}(T_f) = \sigma_{\rm ess}(T_{\tilde f}) \subset \cl \tilde f(\C\setminus D_r)$ for all $r>0$ by the above argument applied to $\tilde f$. For the other inclusion (which involves no product formulas), see~\cite{Stroethoff}. 

For an extension to other Fock spaces 
$$
	F^p_\varphi = \left\{ f\in H(\C) : \int_\C |f(z)|^p e^{-p\varphi(z)}\,dA(z)<\infty\right\}
$$ 
with more general weights $\varphi$ and $0<p<\infty$, see~\cite{HuV-JGEA}, which deals with the so-called doubling weights. These are very general weights that include all standard weights (i.e., $\varphi(z)=-\frac\alpha2 |z|^2$ with $\alpha>0$), the so-called Fock-Sobolev weights, and the weights $\varphi$ for which there are positive constants $m$ and $M$ (depending on $\varphi$) such that
\begin{equation}\label{e:generalized Fock spaces}
	m \le \Delta \varphi \le M
\end{equation}
on $\C$, where $\Delta$ is the Laplacian. It is worth noting that, unlike in these other Fock spaces, we do not currently know whether Fredholmness of Toeplitz operators on doubling Fock spaces can be extended to $\C^n$ due to the lack of suitable estimates for the reproducing kernel.

Let $f\in L^\infty(\C)^{N\times N}$ and consider the block Toeplitz operator $T_f$ on $F^2_N=\{(f_1, \ldots, f_N)^\top : f_j\in F^2\}$. As in the previous two function spaces, the study of Fredholmness of block Toeplitz operators can be reduced to the scalar-valued case using~\eqref{e:6}. However, similarly to $T_f$ on $A^2_N$ with $f\in (C(\overline{\D}) + H^\infty)^{N\times N}$, the index computation in the Fock space setting cannot be reduced to the scalar-valued case and it remains an open problem---perhaps the approach in~\cite{BP2013} can be adapted to this case.

A partial answer to the index computation can be derived from a recent result in~\cite{HV2022}, in which the Schatten class properties of $H_f$ are described in terms of integral distance to analytic functions. More precisely, for $f\in L^2_{\rm loc}(\C)$, define
$$
	G_{r}(f)(z)=\inf_{h\in H(D(z,r))} \left(\frac{1}{|D(z,r)|}\int_{D(z,r)} |f-h|^2 dA \right)^{\frac{1}{2}}\quad (z\in \C),
$$
where $D(z,r)$ is the disk centered at $z$ with radius $r$. For $0<s\le \infty$, we say $f\in \IDA^s$ if $\|G_r(f)\|_{L^s(\C)} < \infty$ for some $r>0$. Notice that the space $\IDA^s$ is independent of $r$. In~\cite{HV2022}, for $0<p<\infty$, it was shown that $H_f$ is in the Schatten class $S_p$ if and only if $f\in \IDA^p$. Let $f\in (L^\infty(\C) \cap \IDA^1)^{N\times N}$ and suppose that $\widetilde \det f$ is bounded away from zero on $\C\setminus \D_R$ for some $R>0$. Then~\eqref{e:6} can be used to show that the entries of $T_f$ commute modulo trace class operators, and hence using the scalar-valued case (see~\cite{BC1987}), we conclude that
$$
	\ind T_f = \ind T_{\det f} = - \wind(\det f|_{|z|=R}).
$$
This result is unsatisfactory because there are bounded symbols that generate compact Hankel operators but do not belong to $\IDA^1$, and further work is required as indicated above.

\section{Sarason's product problem}

In~\cite{Sarason}, Sarason proposed the problem of characterizing the pairs of functions $f, g$ in $H^2$ such that the operator $T_fT_{\bar g}$ is bounded on $H^2$. Related to the present work, he remarked that the identity
\begin{equation}\label{e:Sarason product}
	H^*_{\bar f} H_{\bar g} = T_{f\bar g} - T_fT_{\bar g}
\end{equation}
reduces the problem to the question of when $H^*_{\bar f} H_{\bar g}$ is bounded under the assumption that $fg$ is bounded. When the boundedness assumption on $fg$ is dropped, it can be easily seen that the latter problem is more general (e.g., choose an unbounded $f$ such that $H_{\bar f}$ is bounded and take $g=f$). The following conjecture is often referred to as Sarason's conjecture: \textit{For two functions $f,g$ in $H^2$, $T_fT_{\bar g}$ is bounded if and only if}
\begin{equation}\label{e:Sarason}
	\sup_{z\in\D} \widehat{|f|^2}(z) \widehat{|g|^2}(z) < \infty,
\end{equation}
where $\widehat h$ is defined as the Poisson extension of $h\in L^1(\T)$. In fact, Treil had communicated an argument showing that \eqref{e:Sarason} is necessary to Sarason (see Comment~6 in~\cite{Sarason}) and subsequently Zheng~\cite{Zheng} proved that~\eqref{e:Sarason} with $2$ replaced by $2+\epsilon$ is sufficient. Finally, in the well-known unpublished manuscript of Nazarov~\cite{Nazarov}, it was shown that Sarason's conjecture fails.

A related conjecture was formulated in the Bergman space setting: \textit{For $f,g\in A^2$, $T_fT_{\bar g}$ is bounded on $A^2$ if and only if} 
\begin{equation}\label{e:Sarason2}
	\sup_{z\in\D} \widetilde{|f|^2}(z) \widetilde{|g|^2}(z)<\infty.
\end{equation}
This conjecture was also shown to be false by Aleman, Pott, and Reguera~\cite{APR} using harmonic analysis. However, Stroethoff and Zheng~\cite{SZ2002} showed that if we consider the question of whether $T_f T_{\bar g}$ is both bounded and invertible, then \eqref{e:Sarason} and \eqref{e:Sarason2} provide the right conditions in the settings of $H^2$ and $A^2$, respectively. More precisely, they showed that for $f, g\in A^2$, $T_{f} T_{\bar g}$ is bounded and invertible on $A^2$ if and only~\eqref{e:Sarason2} holds and $\inf \{|f(z)| |g(z)| : z\in \D\}>0$. They also remarked that a similar approach yields an analogous result for Toeplitz operators on the Hardy space, that is, for $f,g\in H^2$, $T_{f} T_{\bar g}$ is bounded and invertible on $H^2$ if and only~\eqref{e:Sarason} holds and $\inf \{|f(z)| |g(z)| : z\in \D\}>0$. It should be noted that the latter result was proved earlier for a pair of outer functions $f,g\in H^2$ by Cruz-Uribe~\cite{CU} using a characterization of invertible Toeplitz operators due to Devinatz and Widom (see, e.g., Theorem~2.23 of~\cite{BS2006}).

Finally, using a number of product identities, Stroethoff and Zheng~\cite{SZ2002} proved that $T_fT_{\bar g}$ is bounded and Fredholm on $A^2$ if and only if~\eqref{e:Sarason2} holds and $\inf_{z\in \D\setminus r\D} |f(z)g(z)| > 0$ for some $r<1$. Again, the same is true in the setting of the Hardy space---just replace~\eqref{e:Sarason2} by~\eqref{e:Sarason}.

Above we have considered Sarason's problem only rather superficially, and while the product formula in~\eqref{e:Sarason product} gives a more general problem involving Hankel operators, the product formulas do not contribute to the two important counterexamples. It is also worth noting that, despite the considerable progress, Sarason's product problem still remains open in the Hardy and Bergman space settings.

We now turn our attention to the Fock space, where Sarason's problem has a simple solution. Indeed, in~\cite{CPZ2014}, for $f,g\in F^2$, it is shown that $T_f T_{\bar g}$ is bounded on $F^2$ if and only if there are $a,b,c\in\C$ such that $f(z) = e^{a+cz}$ and $g(z)=e^{b-cz}$ for all $z\in\C$. One of the key observations is that, when $a\in\C$, $f(z) = e^{\frac12\bar a z}$, $g(z) = e^{-\frac12 \bar a z}$, we have
$$
	T_f T_{\bar g} = e^{\frac14|a|^2} U_a,
$$
where $U_a$ is the unitary operator on $F^2$ defined by
$$
	U_a f(z) = f(z-a) k_a(z)
$$
and $k_a$ is the normalized reproducing kernel of $F^2$ defined by
$$
	k_a(z) = e^{\frac12\bar a z-\frac14|a|^2}.
$$

As weighted Fock spaces $F^2_\varphi$ have received significant attention recently, it would be interesting to know whether something similar holds true for more general weights than those considered in~\cite{BYZ, CPZ2014}. A possible starting point may be the weights $\varphi$ whose Laplacians are bounded above and below (see~\eqref{e:generalized Fock spaces} and~\cite{HuV2022}). What makes Sarason's product problem interesting in this generalized setting is that the reproducing kernel of $F^2_\varphi$ has no explicit representation (unlike in $F^2$) and the unitary operators $U_a$ can no longer be employed. The former obstacle may be possible to overcome with the use of estimates for the (normalized) reproducing kernel, but overall the generalized Sarason's product problem seems nontrivial in generalized Fock spaces and requires new ideas.

\section{Quantization}\label{quantization}

As an application of product formula~\eqref{e:6} and recent work on Hankel operators, we consider deformation quantization (in the sense of Rieffel) and one of its essential ingredients involving the limit condition
\begin{equation}\label{e:Q}
	\lim_{t\to 0} \left \|T^{(t)}_f T^{(t)}_g- T^{(t)}_{fg}\right \|_{{  F}^2_t(\varphi) \to {  F}^2_t(\varphi)} =0,
\end{equation}
where the Toeplitz operators $T^{(t)}_f$ and the Fock spaces $F^2_t(\varphi)$ are defined as follows. For $t>0$, we set
$$
	d\mu_t (z)= \frac 1{t^n} \exp\left \{-  2\varphi\left( \frac    z { \sqrt{t}}\right)\right \} dv(z)
$$
and denote by $L^2_t(\varphi)$ the space of all  Lebesgue measurable functions $f$ in $\C^n$ such that
$$
	\left\|f\right \|_{t} = \left\{ \int_{\C^n} \left|f \right|^2 d\mu_t(z) \right\}^{\frac 12}.
$$
Further, we let $F^2_t(\varphi) = L^2_t(\varphi) \cap H(\C^n)$ and define the Toeplitz operator $T^{(t)}_f$ on $F^2_t(\varphi)$ by
$$
	T^{(t)}_f = P^{(t)} M_f,
$$
where $P^{(t)}$ is the orthogonal projection of $L^2_t(\varphi)$ onto $F^2_t(\varphi)$.

Using the dilation $U_t:  f \mapsto  f(\cdot  \sqrt{t})$, it can be easily shown that
\begin{equation}\label{u-operator-z}
	\| H^{(t)}_f\|_{{  F}^2_t(\varphi) \to {  L}^2_t(\varphi)} 
	= \| H_{f(\cdot \sqrt{t})}\|_{{  F}^2(\varphi) \to {  L}^2(\varphi)},
\end{equation}
where $H^{(t)}_f= (I- P^{(t)})M_f$ is the Hankel operator. To study the limit condition in~\eqref{e:Q}, define for $f\in L^2_{\mathrm{loc}}$, $z\in \C^n$, and $r>0$, 
$$
     MO_{2, r}(f)(z) =  \left( \frac 1{|B(z, r)|} \int_{B(z, r)}
     \left| f-f_{B(z, r)}\right|^2 dv\right)^{\frac 12}
$$
where $B(z,r)= \{ w\in \C^n : |z-w|<r\}$, $f_{S} = \frac 1{|S|} \int_S f dv$ for $S\subset \C^n$ measurable and $dv$ is the usual Lebesgue measure on $\C^n$. Now, let $f\in L^2_{\rm loc}$. We say that $f\in \VMO$ if
 $$
      \lim_{r\to 0} \sup_{z\in \C^n} MO_{2, r}(f)(z)=0.
$$
Further, we say that $f\in \VDA_*$ if 
$$
   \lim_{r\to 0} \sup_{z\in \C^n } G_{2, r}(f)(z)=0.
$$

In~\cite{HuV2022-2}, it was shown that, given $f\in L^\infty$, then for all $g\in L^\infty$, the limit condition in~\eqref{e:Q} holds if and only if $f\in \overline {\mathrm{VDA}}_*$. 

To verify this, notice first that~\eqref{e:6} gives
$$
	T^{(t)}_f T^{(t)}_g- T^{(t)}_{fg} =-\left(H^{(t)}_{\overline f}\right)^* H^{(t)}_g.
$$
for all $f,g \in L^\infty$. Let $f\in \overline{\VDA}_*$. Then, for all $g\in L^\infty$,
\begin{align*}
	\left \|  T^{(t)}_f T^{(t)}_g- T^{(t)}_{fg}\right \|_{{  F}^2_t(\varphi) \to {  F}^2_t(\varphi)} 
	&\le  \| g\|_{L^\infty}\,  \left \| \left(H^{(t)}_{\overline f}\right)^*\right \|_{{  L}^2_t(\varphi)\to {  F}^2_t(\varphi)}\\
	&\le C \|  G_{2, 1}(f(\cdot \sqrt{t}))\|_{L^\infty}\\
	&= C \|G_{2, \sqrt{t}} (f)(\cdot\sqrt{t})\|_{L^\infty}\to 0
\end{align*}
as $t\to 0$, where we used the norm estimate for Hankel operators given in Theorem~1.1 of~\cite{HuV2022-2}. For the converse, again by product formula~\eqref{e:6}, we have
\begin{align*}
	\lim_{t\to 0} \left \|H_{\overline f}^{(t)}\right \|^2_{{  F}^2_t(\varphi) \to {  L}^2_t(\varphi)} 
	&=\lim_{t\to 0} \left \|\left(H_{\overline f}^{(t)} \right)^* H_{\overline f}^{(t)}\right \|_{{  F}^2_t(\varphi) \to {  F}^2_t(\varphi)}\\
	&= \lim_{t\to 0} \|T^{(t)}_f T^{(t)}_{\overline f}- T^{(t)}_{|f|^2}\|_{{  F}^2_t(\varphi) \to {  F}^2_t(\varphi)}  =0,
\end{align*}
and it remains to notice that
$$
	\frac 1 {C}  \|  G_{2, 1}(f(\cdot \sqrt{t}))\|_{L^\infty}
	\le  \left \| \left(H^{(t)}_{\overline f}\right)^*\right \|_{{  L}^2_t(\varphi)\to {  F}^2_t(\varphi)},
$$
which follows from the estimate for Hankel operators mentioned above.

Combining the characterization for~\eqref{e:Q} with the observation that $\VMO = \VDA_* \cap \overline{\VDA}_*$ gives the main result of~\cite{BCH18} (where it was assumed that $\varphi(z)= \frac 1 8 |z|^2$ is the standard weight), that is, given $f\in L^\infty$, then, for all $g\in L^\infty$, it holds that
\begin{equation}\label{T-VMO-a}
        \lim_{t\to 0} \left \|T^{(t)}_f T^{(t)}_g- T^{(t)}_{fg}\right \| =0
        \quad {\rm and}\quad
        \lim_{t\to 0}\left \| T^{(t)}_g T^{(t)}_f - T^{(t)}_{fg}\right \| =0
\end{equation}
if and only if $g\in \mathrm{VMO}$. Here $\|\cdot\| = \|\cdot\|_{{  F}^2_t(\varphi) \to {  F}^2_t(\varphi)}$. For further details, see~\cite{HuV2022-2}.

As for an open problem in this line of work, it would be interesting to characterize those symbols $f\in L^\infty(\C^n)$ for which~\eqref{e:Q} holds for all $g\in L^\infty(\C^n)$ when the operator norm is replaced by the Hilbert-Schmidt (or other Schatten class) norm.

\section*{Acknowledgments}
\noindent I wish to thank Albrecht B\"ottcher for many valuable suggestions. I also thank the American Institute of Mathematics and the Mathematical Sciences Research Institute for their invitations to two meetings in Spring 2022 during which this note was written. This work was supported in part by the American Institute of Mathematics SQuaRE program ``Asymptotic behavior of Toeplitz and Toeplitz+Hankel determinants.''


\begin{thebibliography}{99}
\bibitem{APR} Aleman, A., Pott, S., Reguera, M.C.: Sarason conjecture on the Bergman space. Int. Math. Res. Not. IMRN \textbf{14}, 4320--4349 (2017)

\bibitem{BBE} Basor, E., B\"ottcher, A., Ehrhardt, T.: Harold Widom's contributions to the spectral theory and asymptotics of Toeplitz operators and matrices. This volume.

\bibitem{BC1987} Berger, C.A., Coburn, L.A.: Toeplitz operators on the Segal-Bargmann space. Trans. Amer. Math. Soc. \textbf{301}, 813--829 (1987)

\bibitem{BCH18} Bauer, W., Coburn, L.A., Hagger, R.: Toeplitz quantization on Fock space. J. Funct. Anal. \textbf{274}, 3531--3551 (2018)

\bibitem{BYZ} Bommier-Hato, H., Youssfi, E.H., Zhu, K.: Sarason's Toeplitz product problem for a class of Fock spaces. Bull. Sci. Math. \textbf{141}, 408--442 (2017)

\bibitem{BP2013} Böttcher, A., Perälä, A.: The index formula of Douglas for block Toeplitz operators on the Bergman space of the ball. In: Operator theory, pseudo-differential equations, and mathematical physics, pp. 39--55. Oper. Theory Adv. Appl., 228, Birkhäuser/Springer Basel AG, Basel (2013)

\bibitem{BS2006} B\"ottcher, A., Silbermann, B.: Analysis of Toeplitz operators. Second edition. Springer Monographs in Mathematics. Springer-Verlag, Berlin (2006)

\bibitem{CPZ2014} Cho, H.R., Park, J.-D., Zhu, K.: Products of Toeplitz operators on the Fock space. Proc. Amer. Math. Soc. \textbf{142}, 2483--2489 (2014)

\bibitem{CU} Cruz-Uribe, D.: The invertibility of the product of unbounded Toeplitz operators. Integral Equations Operator Theory \textbf{20}, 231--237 (1994)

\bibitem{HV2022} Hilberdink, T., Virtanen, J.A.: Fredholm theory of Toeplitz operators on $H^1$. II. (in preparation)

\bibitem{HuV2022} Hu, Z., Virtanen, J.A.: Schatten class Hankel operators on the Segal-Bargmann space and the Berger-Coburn phenomenon. Trans. Amer. Math. Soc. \textbf{375}, 3733-3753 (2022)

\bibitem{HuV-JGEA} Hu, Z., Virtanen, J.A.: Fredholm Toeplitz operators on doubling Fock spaces. J. Geom. Anal. \textbf{32}, 106 (2022)

\bibitem{HuV2022-2} Hu, Z., Virtanen, J.A.: IDA and Hankel operators on Fock spaces. Anal. PDE (in press) arXiv:2111.04821

\bibitem{McDonald} McDonald, G.: Fredholm properties of a class of Toeplitz operators on the ball. Indiana Univ. Math. J. \textbf{26}, 567--576 (1977)

\bibitem{Nazarov} Nazarov, F.: A counterexample to Sarason’s conjecture. Unpublished manuscript (1997) \url{https://users.math.msu.edu/users/fedja/prepr.html}

\bibitem{Sarason} Sarason, D.: Products of Toeplitz operators. In: V. P. Khavin, N. K. Nikolski (eds.), Linear and Complex Analysis Problem Book 3, Part I, Lecture Notes in Math., 1573, pp. 318--319. Springer, Berlin (1994)

\bibitem{Stroethoff} Stroethoff, K.: Hankel and Toeplitz operators on the Fock space. Michigan Math. J. \textbf{39}, 3--16 (1992)

\bibitem{SZ2002} Stroethoff, K., Zheng, D.: Invertible Toeplitz products. J. Funct. Anal. \textbf{195}, 48--70 (2002)

\bibitem{V2006} Virtanen, J.A.: Fredholm theory of Toeplitz operators on the Hardy space $H^1$. Bull. London Math. Soc. \textbf{38}, 143--155 (2006)

\bibitem{Widom} Widom, H.: Asymptotic behavior of block Toeplitz matrices and determinants. II. Advances in Math. \textbf{21}, 1--29 (1976)

\bibitem{Zheng} Zheng, D.: The distribution function inequality and products of Toeplitz operators and Hankel operators. J. Funct. Anal. \textbf{138}, 477--501 (1996)
\end{thebibliography}
\end{document}